\documentclass[11pt]{article}
\usepackage[ansinew]{inputenc}
\usepackage{array}
\usepackage{caption}
\usepackage{color}
\usepackage{amsmath}
\usepackage{amsthm}
\usepackage{titlesec}
\usepackage{amsxtra}
\usepackage{amstext}
\usepackage{amssymb}
\usepackage{latexsym}
\usepackage{float}
\usepackage{graphicx}
\usepackage{amsfonts,cite}
\usepackage{graphics,epstopdf}
\usepackage{epsfig}
\usepackage[colorlinks=true,allcolors=blue]{hyperref}
\usepackage{cite}
\theoremstyle{definition}

\numberwithin{equation}{section}
	
\titleformat{\section}[hang]{\Large\bfseries}{\thesection.}{1em}{}

\begin{document}

\begin{center}
{\huge \textbf{ Some properties of general-$\lambda$-matrix polynomials: an umbral approach}}\\~

{\bf Ghazala Yasmin*${}^{}${\footnote{$^{}$*Corresponding author; E-mail:~ghazala30@gmail.com (Ghazala Yasmin)}} and  Aditi Sharma${}^{}$ {\footnote{$^{}$ E-mail:~aditisharma@gmail.com( Aditi Sharma)}}    }\\
Department of Applied Mathematics, Aligarh Muslim University, Aligarh, India \\
\end{center}
\parindent=8mm
\noindent
\begin{abstract} \noindent
``2-variable general-$\lambda$-matrix polynomials (2VG$\lambda$MP)" is a new family of matrix polynomials, introduced and studied in this article. These matrix polynomials are constructed using umbral and symbolic methods. We delve into the generating function, explicit series representation, differential equation, quasi-monomiality, summation formula, determinant representation and various identities satisfied by these polynomials. Furthermore, a thorough investigation are conducted to ascertain the outcomes for the members of the general-$\lambda$-matrix family. Additionally, 3D graphs and figures showing the zeros distribution, real zeros and stacks of zeros for few members of the family of these polynomials are also presented in the article.
\end{abstract}
\vspace{-4pt}\noindent
{\bf{Mathematics Subject Classifications:}}~~ 33C45, 33C47, 33F10, 33D05.\\

\vspace{-9pt}\noindent
{\bf{\em Keywords:}} 2-variable general polynomials, $\lambda$-matrix polynomials, Symbolic operators, Umbral methods, Truncated exponential-$\lambda$-matrix polynomials.
\vspace{-12pt}
\section{Introduction and Preliminaries}
\vspace{-7pt}
Umbral calculus is a powerful tool in the study of special functions, often encountered in mathematical physics and applied mathematics, such as Bessel functions, Legendre polynomials and hypergeometric functions. Using symbolic methods, umbral calculus abstracts the properties and operations of these functions, simplifying the process of proving identities, deriving recurrence relations and exploring integral representations or asymptotic behaviors. This approach, which emphasizes algebraic and combinatorial aspects through symbolic manipulation, enables efficient management of complex expressions and often reveals new insights that might not be apparent through traditional methods.

\vspace{+4pt} \noindent
One of the notable strengths of umbral calculus is its effectiveness in analyzing the generating functions, which are commonly used to represent and study sequences related to special functions. By manipulating these generating functions symbolically, novel results and uncovering interrelationships among various special functions can be established. Umbral calculus also facilitates the exploration of generalized functions and provides a unified framework for examining series expansions. This methodology not only streamlines complex mathematical tasks but also enhances the ability to detect underlying patterns and relationships within special functions. Ultimately, umbral calculus simplifies proofs and derives new results with greater ease compared to conventional approaches, broadening the scope and depth of mathematical analysis. Symbolic methods, in general, involve manipulating abstract symbols rather than specific numerical values, allowing for generalization and simplification by focusing on the relationships and operations between symbols. This approach helps in revealing the intricate structures and patterns inherent in mathematical problems.

\vspace{+4pt} \noindent
Recently, Dattoli et. al \cite{dattoli2017circular} introduced a new family of polynomials, known as ``the $\lambda$-polynomials $\lambda_n(x,y)$", which have the potential to significantly contribute to the understanding and analysis of special functions. The following generating function is used to define the $\lambda$-polynomials $\lambda_n(x,y)$\cite{dattoli2017circular}:
\begin{equation}
e^{yt}cos{\sqrt{xt}}=\sum_{n=0}^{\infty}\lambda_n(x,y)\frac{t^n}{n!},
\end{equation}
and possess the following series expansion:
\begin{equation}
  \lambda_n(x,y)=n!\sum_{j=0}^{n}\frac{(-1)^jx^jy^{n-j}}{(2j)!(n-j)!}.  
\end{equation}
Matrix polynomials are important in extending scalar polynomial concepts to matrix-valued functions, broadening the scope of special functions. The expansion of orthogonal polynomials to matrix settings has diverse applications in quantum mechanics and control theory. Matrix polynomials also play a key role in special functions, including matrix hypergeometric and Bessel functions and are essential for solving differential equations and analyzing symmetries. In combinatorics and probability theory, they are crucial for generating functions and stochastic processes, providing a versatile framework across various mathematical and applied fields.

\vspace{+4pt} \noindent
The extension and generalization of $\lambda$-polynomials in which matrix is used as a parameter have been recently established and investigated for matrices in ${\mathbb{C}}^{\mathbb{N}*\mathbb{N}}$.
The following is the definition for $\lambda$-matrix polynomials in symbolic terms\cite{zainab2024symbolic}:
\begin{equation}
\lambda_n^{R}(x,y)=\hat{J}^{R}(y-\hat{J}x)^n\psi_0\hspace{0.08cm},\hspace{1cm}{x,y}\in \mathbb{R} , \forall n \in \mathbb{N},
\end{equation}
where Dattoli et al.\cite{dattoli2017circular} provide a symbolic operator $\hat{J}$ that acts on the vacuum function $\psi_z$ as:
\begin{equation}
\hat{J}^j\psi_z={\Gamma(j+z+1)}{[\Gamma(2(j+z)+1)]^{-1}}, \hspace{1cm} j\in \mathbb{R} ,
\end{equation}
such that
\begin{equation}
\hat{J}^p\hat{J}^q=\hat{J}^{p+q}.
\end{equation}
The ordinary and exponential generating functions and the explicit series representation of $\lambda$-matrix polynomials in two variables are defined as\cite{zainab2024symbolic}:
\begin{equation}
   \sum_{n=0}^{\infty} t^n \lambda_n^{R}(x,y) =\frac{1}{1-yt}e_0\left(\frac{xt}{1-yt};R\right),
\end{equation}
\begin{equation}
 \sum_{n=0}^{\infty} \lambda_n^{R}(x,y)\displaystyle \frac{t^n}{n!} =e^{yt}\cos(\sqrt{xt};R)
 \end{equation}
 and 
 \begin{equation}
 \lambda_n^{R}(x,y) =n!\sum_{j=0}^{n} \displaystyle \frac{(-1)^j\Gamma{(R+(j+1)I)x^jy^{n-j}}}{j!{(n-j)!}\Gamma{(2R+(2j+1)I)}}.
 \end{equation}
where the associated cosine function $\cos(x;R)$; $R\in\mathbb{N}$ is symbolically defined as \cite{zainab2024symbolic}:
\begin{equation}
\cos(x;R)=\hat{J}^{R}e^{-\hat{J}x^2}\psi_0=\sum_{j=0}^{\infty} \displaystyle \frac{(-1)^j\Gamma{(R+(j+1)I)}x^{2j}}{j!\Gamma{(2R+(2j+1)I)}}.
\end{equation}
The $\lambda$-matrix polynomials in one variable $\lambda_n^{R}(x)$ are specified using the following definition \cite{zainab2024symbolic}:
\begin{equation}
\lambda_n^{R}(x)=\hat{J}^R(x-\hat{J})^n\psi_0.\end{equation}
The generating function and explicit series representation of $\lambda_n^R(x)$ is given by \cite{zainab2024symbolic}:
\begin{equation}
\sum_{n=0}^{\infty} \lambda_n^{R}(x)\displaystyle \frac{t^n}{n!}=e^{xt}\cos(\sqrt{t};R)\end{equation}
and
\begin{equation}
\lambda_n^{R}(x) =n!\sum_{j=0}^{n} \displaystyle \frac{(-1)^j x^{n-j}\Gamma{(R+(j+1)I)}}{j!{(n-j)!}\Gamma{(2R+(2j+1)I)}}.
\end{equation}
From equations (1.3) and (1.10), we conclude that 
\begin{equation}
    \lambda_n^{R}(x)=\lambda_n^{R}(1,x).
\end{equation}
\noindent
The condition that $R$ satisfies is as follows: $Re(w)>0$ $\forall$ w$\in\sigma(R)$ , where $R$ is a positive stable matrix in ${\mathbb{C}}^{\mathbb{N}*\mathbb{N}}$ throughout this paper.
When the complex plane $F_0$ cut along the negative real axis then $w^{\frac{1}{2}}$ represents $exp(\frac{1}{2}log(w))$ , where $log(w)$ denotes the principle logarithm of w, . If R is a matrix in $\mathbb{C^{\mathbb{N}\times \mathbb{N}}}$ with $\sigma(R) \subset F_0$, then $R^{\frac{1}{2}}=\sqrt{R}$ denotes the image by $w^{\frac{1}{2}}$ of the matrix functional calculus\cite{dunford1988linear} acting on the matrix R.  

\vspace{+4pt} \noindent
Novel approaches to the analysis of solutions to different PDEs found in physical problems have been made possible by the development of the general class of polynomials of two variables. These polynomials have proven to be instrumental in addressing challenges in optics and quantum mechanics. In view of their broad applications, a specific category known as the 2-variable general polynomials (2VGP) $p_n(x,y)$ \cite{khan2013general} has emerged as a key focus for further exploration and study. These polynomials have the following generating function defined \cite{khan2013general}:
\begin{equation}
 e^{xt} \phi(y,t)= \sum_{n=0}^{\infty} p_n(x,y) \frac{t^n}{n!}, \hspace{1cm} p_0(x,y)=1,
\end{equation}
where the following formal series expansion of $\phi(y,t)$ exists:
\begin{equation}
    \phi(y,t)= \sum_{n=0}^{\infty} \phi_n(y) \frac{t^n}{n!} .
\end{equation} 
Next, a table is given below showcasing ``The family of 2-variable general polynomials" for a comprehensive understanding of their structure and properties.
\pagebreak

\begin{table}[h]

\vspace{.1cm}

\noindent{	
	{\footnotesize{
			\begin{tabular*}{\textwidth}{ @{\extracolsep{\fill}} l@{\hspace{.1cm}}llll}
				\hline
				\bf{S.No.}&{\bf Name of the polynomial} & \bf{Series definition  }&\bf{Generating function }\\
				\hline
				\\
				\bf{I.}&Gould Hopper &$H_n^{m}(x,y)=n!\sum_{j=0}^{[\frac{n}{m}]}\frac{y^j}{j!}\frac{x^{n-mj}}{(n-mj)!}$&$e^{xt} e^{yt^m}=\sum_{n=0}^{\infty} H_n^{(m)}(x,y) \frac{t^n}{n!}$\\
				&polynomials $H_n^{m}(x,y)$\cite{gould1962operational}\\
				\\
				\bf{II.} &Hermite polynomials $H_n(x,y)$\cite{appell1926fonctions}   & $H_n(x,y)=n!\sum_{j=0}^{[\frac{n}{2}]}\frac{y^j}{j!}\frac{x^{n-2j}}{(n-2j)!}$&$e^{xt} e^{yt^2}=\sum_{n=0}^{\infty} H_n(x,y) \frac{t^n}{n!}$ \\
				\\
				\bf{III.} &Laguerre polynomials&$L_n(y,x)=n!\sum_{j=0}^{n}{(-1)^j}\frac{y^j}{(j!)^2}\frac{x^{n-j}}{(n-j)!}$&$e^{xt} C_0(yt)=\sum_{n=0}^{\infty} L_n(y,x) \frac{t^n}{n!}$ \\
				& $L_n(y,x)$\cite{dattoli1999hermite}\cite{dattoli1998operational}\\
				\\
				\bf{IV.} &Generalized Laguerre&$_{m}L_n(y,x)=n!\sum_{j=0}^{[\frac{n}{m}]}\frac{y^j}{(j!)^2}\frac{x^{n-mj}}{(n-mj)!}$&$ e^{xt} C_0(-yt^m)=\sum_{n=0}^{\infty} {_{m}L_n(y,x)} \frac{t^n}{n!}$\\
				&polynomials $_{m}L_n(y,x)$\cite{dattoli1999generalized}\\
				\\
				\bf{V.} &Truncated exponential&$e_n(x,y)=n!\sum_{j=0}^{n}\frac{x^{n-j}~y^j}{(n-j)!}$&$\frac{e^{xt}}{(1-yt)}=\sum_{n=0}^{\infty}e_n(x,y)\frac{t^n}{n!}$\\
				&polynomials of order 1&&\\
				& $e_n(x,y)$\cite{dattoli2004class} \\
				\\
		\hline
\end{tabular*}}}}\\

\caption{\textbf{Family of the 2-variable general polynomials}}

\end{table}

\vspace{+4pt} \noindent
In the 1940s, J.F. Steffensen\cite{steffensen1941poweroid} developed the concept of monomiality through his work on the poweroid. Dattoli redeveloped the monomiality principle\cite{dattoli1999hermite}, which asserts that two operators, a multiplicative operator $\hat{M}$ and a derivative operator $\hat{P}$, are necessary for a polynomial set $Y_n(y)$ to be considered a ``quasi monomial" if they satisfy the following recurrence relations:
\begin{equation}
\hat{M}{(Y_n(y))}=Y_{n+1}(y)
\end{equation}
and
\begin{equation}
\hat{P}{(Y_n(y))}=nY_{n-1}(y).
\end{equation}
The following commutation property is met by the operators $\hat{M}$ and $\hat{P}$:
 \begin{equation}
 [\hat{P},\hat{M}]=\hat{P}\hat{M}-\hat{M}\hat{P}=\hat{1},
 \end{equation}
and fulfill the structure of a Weyl group.

\vspace{+4pt} \noindent
Moreover, the polynomial $Y_n(y)$ satisfies the following differential equation if $\hat{P}$ and $\hat{M}$ have a differential realization:
\begin{equation}
\hat{M}\hat{P}(Y_n(y))=nY_n(y).
\end{equation}
Motivated and inspired by the work mentioned above, we introduce and explore a new family of general $\lambda$-matrix polynomials and establish certain interesting results.
Umbral methods are used in section 2 to introduce the 2-variable general-$\lambda$-matrix polynomials. The properties of these polynomials, including their generating function, series definition, summation formulas, integral representation, differential recurrence relations and the numerous identities they satisfy, are covered in great detail. Section 3 delve into the discussion of quasi-monomiality and determinant form of 2VG$\lambda$MP $_{G}\lambda_n^R(x,y)$. Finally, section 4 of this article gives a detailed analysis of several examples of 2VG$\lambda$MP $_{G}\lambda_n^R(x,y)$ which also includes 3D graphs and zeros distribution of these polynomials.
\pagebreak

\vspace{-12pt}
\section{\texorpdfstring{2-Variable General-$\lambda$-Matrix Polynomials (2VG$\lambda$MP) $_{G}\lambda_n^R(x,y)$}{2-Variable General-lambda-Matrix Polynomials (2VGlambdaMP) G-lambda-n-R(x,y)}}

\vspace{-7pt}
This section starts by defining the 2VG$\lambda$MP $_{G}\lambda_n^R(x,y)$ as the discrete $\lambda$-matrix convolution of the 2VG$\lambda$MP $p_n(x,y)$ and establishes certain properties for these polynomials. To achieve this aim we first define a symbolic operator $\hat{q_y}$ which acts on vacuum $\xi_0$ such that \begin{equation}\phi_n(y)=\hat{q_y^n}\xi_0.\end{equation}
Thus,
\begin{equation}
\phi(y,t)=\sum_{n=0}^{\infty}\hat{q_y^n}\xi_0\displaystyle\frac{t^n}{n!}=e^{(\hat{q_yt})}\xi_0,
\end{equation}
which on using $\phi(y,t)=e^{\hat{q_yt}}\xi_0 $ in eqn (1.1) gives the symbolic definition of 2-variable general polynomials $p_n(x,y)$ as:
\begin{equation}
 p_n(x,y)=(x+\hat{q_y})^n\xi_0.\end{equation}
 Next, we introduce the 2-variable general-$\lambda$-matrix polynomials (2VG$\lambda$MP) $_{G}\lambda_n^R(x,y)$ as:
 \begin{equation}
     {_{G}\lambda_n^R(x,y)}={_{G}\lambda_n^R(x+\hat{q}_y)\xi_0}.
 \end{equation}
 On using equation (1.10) in (2.4), we get symbolic definition of 2VG$\lambda$MP $_{G}\lambda_n^R(x,y)$ as:
\begin{equation}
_{G}\lambda_n^{R}(x,y)=\hat{J}^R(x+\hat{q}_y-\hat{J})^n\xi_0\psi_0,
\end{equation}
 where, $\hat{J}$ operates on $\psi_0$ and $\hat{q}_y$ operates on $\xi_0$.

\vspace{+4pt} \noindent
Now we examine some of the properties which the 2VG$\lambda$MP satisfies. It is established by proving the following result that the 2VG$\lambda$MP exhibits a generating relation.

\vspace{+4pt} \noindent
\textbf{Theorem 2.1.} The generating relations for the 2-variable general-$\lambda$-matrix polynomials 2VG$\lambda$MP ${_{G}\lambda_n^{R}(x,y)}$ are as follows:
\begin{equation} 
\sum_{n=0}^{\infty}{_{G}\lambda_n^{R}(x,y)} \displaystyle\frac{t^n}{n!}=e^{xt}\phi(y,t)\cos({\sqrt{t};R}).
\end{equation}

\vspace{+4pt} \noindent
\textbf{Proof.}
From equation (2.5), we have
\begin{equation}
 \sum_{n=0}^{\infty}{_{G}\lambda_n^{R}(x,y)} \displaystyle\frac{t^n}{n!}= \hat{J}^Re^{(x+\hat{q}_y-\hat{J})t}\hspace{0.05cm}\xi_0\psi_0.
 \end{equation}
In view of Weyl's identity\cite{louisell1973quantum}, the above equation leads to
\begin{equation}
 \sum_{n=0}^{\infty}{_{G}\lambda_n^{R}(x,y)} \displaystyle\frac{t^n}{n!}= \hat{J}^Re^{xt}e^{\hat{q}_yt}e^{-\hat{J}t}\hspace{0.05cm}\xi_0\psi_0,
 \end{equation}
which yields assertion (2.6), on using definitions (1.9) and (2.2) in the previous equation.

\vspace{+4pt} \noindent 
 \textbf{Remark 1.} We observe from equation (2.6) that:
 \begin{equation}
 t{_{G}\lambda_n^{R}(x,y)}=n{_{G}\lambda_{n-1}^{R}(x,y)}.
 \end{equation}
Next, by demonstrating the following result, the explicit series representation of the 2VG$\lambda$MP$ {_{G}\lambda_n^{R}(x,y)}$ is obtained:

\vspace{+4pt} \noindent
 \textbf{Theorem 2.2.} The following series expansion is valid for the 2VG$\lambda$MP ${_{G}\lambda_n^{R}(x,y)}$:
 \begin{equation} {_{G}\lambda_n^{R}(x,y)}=\sum_{j=0}^{n} \binom{n}{j} (-1)^jp_{n-j}(x,y)\Gamma{(R+(j+1)I)}(\Gamma{(2R+(2j+1)I))^{-1}}.
\end{equation}
\vspace{+4pt} \noindent
\textbf{Proof.}
Extending the r.h.s of equation (2.5), it is evident that
\begin{equation} {_{G}\lambda_n^{R}(x,y)}=\hat{J}^R\sum_{j=0}^{n}\binom{n}{j} {(x+\hat{q}_y)}^{n-j}(-\hat{J})^j\xi_0\psi_0.
\end{equation}
Assertion (2.10) is made using symbolic definitions (2.3), (1.4) and (1.5) in equation (2.11).

\vspace{+4pt} \noindent
We demonstrate the following outcome to determine the differential recurrence relation of the 2VG$\lambda$MP ${_{G}\lambda_n^{R}(x,y)}$:

\vspace{+4pt} \noindent
\textbf{Theorem 2.3.} The following differential recurrence relation is valid for the 2VG$\lambda$MP ${_{G}\lambda_n^{R}(x,y)}$:
\begin{equation}
\displaystyle\frac{\partial}{\partial x}{_{G}\lambda_n^{R}(x,y)}=n{_{G}\lambda_{n-1}^{R}(x,y)}.
\end{equation}
More generally,
\begin{equation}
\displaystyle\frac{\partial^j}{{\partial x}^j}{_{G}\lambda_n^{R}(x,y)}=\displaystyle\frac{n!}{(n-j)!} {_{G}\lambda_{n-j}^{R}(x,y)}.
\end{equation}
\vspace{+4pt} \noindent
\textbf{Proof.} Differentiating  (2.6), w.r.t $x$, we get
\begin{equation} \frac{\partial}{\partial x}\sum_{n=0}^{\infty}{_{G}\lambda_n^{R}(x,y)} \displaystyle\frac{t^n}{n!}=te^{xt}\phi(y,t)\cos({\sqrt{t};R}).
\end{equation}
Applying the generating function (2.6) to the r.h.s. of the previous equation and simplifying yields
\begin{equation} \frac{\partial}{\partial x}\sum_{n=0}^{\infty}{_{G}\lambda_n^{R}(x,y)} \displaystyle\frac{t^n}{n!}= n\sum_{n=1}^{\infty}{_{G}\lambda_{n-1}^{R}(x,y)} \displaystyle\frac{t^{n}}{n!}.
\end{equation}
Coefficients of similar powers of t on both sides of the expression (2.15) are compared to arrive at the claim (2.12).
Moreover, mathematical induction is applied for the proof of (2.13). It is evident from the identity (2.12) that the result (2.13) is true for $j=1$.
From the principle of mathematical induction, assertion (2.12) follows assertion (2.13).

\vspace{+4pt} \noindent
\textbf{Theorem 2.4.} The explicit representation of the 2VG$\lambda$MP ${_{G}\lambda_n^{R}(x,y)}$ is as follows:
\begin{equation}
    {_{G}\lambda_n^{R}(x,y)}=\sum_{j=0}^{n}\binom{n}{j}{\lambda_j^{R}(x)}{\phi_{n-j}(y)}.
\end{equation}
\vspace{+4pt} \noindent
\textbf{Proof.} Make use of the generating relations given by equations (1.11), (1.15) and (2.6) to obtain
\begin{equation}
\sum_{n=0}^{\infty}{_{G}\lambda_n^{R}(x,y)}\frac{t^n}{n!}=\sum_{j=0}^{\infty}{\lambda_j^{R}(x)}\frac{t^j}{j!}\sum_{n=0}^{\infty}\phi_n(y)\frac{t^n}{n!},
\end{equation}
which on using Cauchy product rule on r.h.s. gives
 \begin{equation}
  \sum_{n=0}^{\infty}{_{G}\lambda_n^{R}(x,y)} \displaystyle\frac{t^n}{n!}= \sum_{n=0}^{\infty}\sum_{j=0}^{n} \binom{n}{j}{\lambda_{j}^{R}(x)}{\phi_{n-j}(y)} \displaystyle\frac{t^n}{n!}.
\end{equation} 
By equating the coefficients of $t$ on both sides in the above equation, assertion (2.16) is obtained.

\vspace{+4pt} \noindent
We now prove the following result to derive the summation formulae for 2-variable general-$\lambda$-matrix polynomials ${_{G}\lambda_n^{R}(x,y)}$.

\vspace{+4pt} \noindent
\textbf{Theorem 2.5.} The following implicit summation formulae hold for the 2-variable general-$\lambda$-matrix polynomials ${_{G}\lambda_n^{R}(x,y)}$:
\begin{equation}
_G\lambda_n^{R}(x+z,y)=\sum_{j=0}^{n}\binom{n}{j}z^j{_G\lambda_{n-j}^{R}(x,y)}.
\end{equation}
\vspace{+4pt} \noindent
\textbf{Proof.} In generating relation (2.7), substituting $x+z$ for $x$, gives
 \begin{equation}
 \sum_{n=0}^{\infty}{_{G}\lambda_n^{R}(x+z,y)} \displaystyle\frac{t^n}{n!}= \hat{J}^Re^{(x+\hat{q}_y-\hat{J})t}e^{zt}\hspace{0.05cm}\xi_0\psi_0,
 \end{equation}
which, when relation (2.7) and the exponential series expansion is applied in the above expression, we find,
 \begin{equation}
 \sum_{n=0}^{\infty}{_{G}\lambda_n^{R}(x+z,y)} \displaystyle\frac{t^n}{n!}=
 \sum_{n=0}^{\infty}\sum_{j=0}^{\infty}{_{G}\lambda_n^{R}(x,y)} \displaystyle\frac{t^n}{n!}\frac{{z^j}{t^j}}{j!}.
 \end{equation}
Furthermore, using the Cauchy product rule in the r.h.s. of equation (2.21) and then comparing the coefficients of identical powers of $t$, assertion (2.19) is proven.

\vspace{+4pt} \noindent
\textbf{Theorem 2.6.} The following identity is true for the 2-variable general-$\lambda$-matrix polynomials ${_{G}\lambda_n^{R}(x,y)}$:
\begin{equation}
_{G}\lambda_{r+s}^{R}(z,y)=\sum_{n,j=0}^{r,s}\binom{r}{n}\binom{s}{j}(z-x)^{n+j}{_{G}\lambda_{r+s-n-j}^{R}(x,y)}.
\end{equation}
\vspace{+4pt} \noindent
\textbf{Proof.} Substituting $t+u$ for $t$ in equation (2.7) and applying the subsequent identity
\begin{equation}
\sum_{m=0}^{\infty}f(m)\frac{(t+u)^m}{m!}=\sum_{m,j=0}^{\infty}f(m+j)\frac{t^m}{m!}\frac{u^j}{j!},
\end{equation}
in the r.h.s. of the resulting expression gives
\begin{equation}
e^{(\hat{q}_y-\hat{J})(t+u)}\hat{J}^R \xi_0 \psi_0 =e^{-x(t+u)} \sum_{r,s=0}^{\infty} {_G\lambda_{r+s}^R}(x,y)\frac{t^r}{r!} \frac{u^s}{s!}.
\end{equation}
On multiplying $e^{z(t+u)}$ on both sides and using identity (2.23) in the resultant equation, we obtain
\begin{equation}
\sum_{r,s=0}^{\infty} {_G\lambda_{r+s}^R}(z,y)\frac{t^r}{r!} \frac{u^s}{s!}=\sum_{n,j=0}^{\infty}\frac{(z-x)^{n+j}t^nu^j}{n!j!} \sum_{r,s=0}^{\infty} {_G\lambda_{r+s}^R}(x,y)\frac{t^r}{r!} \frac{u^s}{s!}.
\end{equation}
Also, by simplifying the resulting equation and applying the Cauchy product rule to the r.h.s. of equation (2.25), we have
\begin{equation}
\sum_{r,s=0}^{\infty} \frac{u^s t^r}{s! r!} {_G\lambda^{R}_{r+s}}(z,y) = \sum_{r,s=0}^{\infty} \sum_{n,j=0}^{r,s} \frac{(z-x)^{n+j} t^r u^s}{n! j! (r-n)! (s-j)!} {_G\lambda^{R}_{r+s-n-j}}(x,y).
\end{equation}
which yields assertion (2.22) when the similar powers of $t$ and $s$ are compared.

\vspace{+4pt} \noindent
\textbf{Theorem 2.7.} The summation formula for the 2-variable general-$\lambda$-matrix polynomial ${_{G}\lambda_n^{R}(x,y)}$ is as follows:
\begin{equation}
{_{G}\lambda_{n+1}^{R}(x+z,y)}-{_{G}\lambda_{n+1}^{R}(x,y)}=\sum_{j=1}^{n}\binom{n+1}{j}{_{G}\lambda^{R}_{n-j+1}(x,y)z^j}.
\end{equation}
\vspace{+4pt} \noindent
\textbf{Proof.} Using equation (2.6), we obtain
\begin{equation}
\sum_{n=0}^{\infty}{_{G}\lambda_n^{R}(x+z,y)}\frac{t^n}{n!}-\sum_{n=0}^{\infty}{_{G}\lambda_n^{R}(x,y)}\frac{t^n}{n!}=(e^{zt}-1)e^{xt}\phi(y,t)\cos(\sqrt{t};R).
\end{equation}
In view of equation (2.6) and series expansion of an exponential function, it gives
\begin{equation}
\sum_{n=0}^{\infty}{_{G}\lambda_n^{R}(x+z,y)}\frac{t^n}{n!}-\sum_{n=0}^{\infty}{_{G}\lambda_n^{R}(x,y)}\frac{t^n}{n!}=\sum_{n=0}^{\infty}{_{G}\lambda_n^{R}(x,y)}\frac{t^n}{n!}\sum_{j=0}^{\infty}\frac{(zt)^{j+1}}{{(j+1)}!}.
\end{equation}
Now, by applying the Cauchy product rule to the r.h.s. of the resulting equation and simplifying, we obtain
 \begin{equation}
\sum_{n=0}^{\infty}{_{G}\lambda_n^{R}(x+z,y)}\frac{t^n}{n!}-\sum_{n=0}^{\infty}{_{G}\lambda_n^{R}(x,y)}\frac{t^n}{n!}= \sum_{n=0}^{\infty}\sum_{j=1}^{n} \binom{n+1}{j}z^j{_{G}\lambda^{R}_{n-j+1}(x,y)} \displaystyle\frac{t^{n+1}}{{(n+1)}!}.
\end{equation} 
As we solve for the coefficients of comparable powers of $t$, we arrive at assertion (2.27).

\vspace{+4pt} \noindent
\textbf{Theorem 2.8.} The integral representation that follows is valid for the 2VG$\lambda$MP:
\begin{equation}
\int_{x}^{x+z} {_G\lambda_n^R}(x,y)dx=\frac{1}{n+1}\sum_{j=1}^{n}\binom{n+1}{j} {_G\lambda_{n-j+1}^R(x,y)}z^j.
\end{equation}
\vspace{+4pt} \noindent
\textbf{Proof.} By applying expression (2.12), we obtain
\begin{align*}
    \int_{x}^{x+z} {_G\lambda_n^R}(x,y)dx&=\frac{1}{n+1}[{_G\lambda_{n+1}^R(x,y)}]_{x}^{x+z}
    \\
    &=\frac{1}{n+1}[{_G\lambda_{n+1}^R(x+z,y)}-{_G\lambda_{n+1}^R(x,y)}],
\end{align*}
which on using equation (2.27) leads us to assertion (2.31).

\vspace{+4pt} \noindent
The quasi-monomiality and determinant form of 2-variable general-$\lambda$-matrix polynomials are discussed in the next section.
\vspace{-12pt} 
\section{Monomiality principle and Determinant form}
\vspace{-7pt} 
In this section, we establish the quasi-monomial characteristics of the 2VG$\lambda$MP ${_{G}\lambda_n^{R}(x,y)}$ which are essential as they advance to the understanding of an umbral aspect of these newly introduced family within the context of monomiality. The validation of these properties enriches the mathematical analysis and makes these polynomials more applicable in various mathematical contexts. In order to prove that the 2-variable general-$\lambda$-matrix polynomials 2VG$\lambda$MP ${_{G}\lambda_n^{R}(x,y)}$ are quasi-monomial, we do the following:

\vspace{+4pt} \noindent
\textbf{Theorem 3.1.} For the 2VG$\lambda$MP, the following multiplicative and derivative operators are obtained.
\begin{equation}
\hat{M}_{(G^{\lambda})}=(x+\hat{q}_y-\hat{J})
\end{equation}
and
\begin{equation}
\hat{P}_{(G^{\lambda})}=D_x,
 \end{equation} respectively. 
 
 \vspace{+4pt} \noindent
 \textbf{Proof.} Operating $(x+\hat{q}_y-\hat{J})$ on both sides of equation (2.5), we obtain
 \begin{equation}
 (x+\hat{q}_y-\hat{J}){_{G}\lambda_n^{R}(x,y)}=\hat{J}^R(x+\hat{q}_y-\hat{J})^{n+1}\xi_0\psi_0,
\end{equation}
which again using equation (2.5), gives
 \begin{equation}
 (x+\hat{q}_y-\hat{J}){_{G}\lambda_n^{R}(x,y)}={_{G}\lambda_{n+1}^{R}(x,y)}.
\end{equation}
Additionally, we are led to assertion (3.1) in view of expressions (1.16) and (3.4). Now, expression (2.5) is differentiated with respect to $x$ and then using equation (2.5) on the r.h.s. of the resulting equation gives
\begin{equation}
D_x{_{G}\lambda_n^{R}(x,y)}=n{_{G}\lambda_{n-1}^{R}(x,y)}.
\end{equation}
Thus in view of equations (1.17) and (3.5), statement (3.2) is derived.

\vspace{+4pt} \noindent
\textbf{Remark 2.} The following result is obtained by applying equations (3.1) and (3.2) to equation (1.19).

\vspace{+4pt} \noindent
 \textbf{Corollary.} The following differential equation is satisfied by the 2-variable general-$\lambda$-matrix polynomials ${_{G}\lambda_n^{R}(x,y)}$:
 \begin{equation}
 ((x+\hat{q}_y-\hat{J})\displaystyle\frac{\partial}{\partial x}-n){_{G}\lambda_n^{R}(x,y)}=0.
 \end{equation}
 The 2-variable general-$\lambda$-matrix polynomials 2VG$\lambda$MP ${_{G}\lambda_n^{R}(x,y)}$ are framed within the symbolic approach. Further, establishing the determinant form of these polynomials is a novel approach to derive the properties employing basic linear algebra tools and has wide applications in diverse areas. We derive the determinant definition of 2VG$\lambda$MP by proving the following results:

\vspace{+4pt} \noindent
\textbf{Theorem 3.2.} The 2VG$\lambda$MP  is defined by
\begin{equation}
{_{G}\lambda_0^{R}(x,y)}=\displaystyle\frac{1}{\gamma_0(y)}\frac{\Gamma(R+I)}{\Gamma(2R+I)}.
\end{equation}
\begin{equation}
_G\lambda_n^{R}(x,y) = \frac{(-1)^n}{\gamma_0{(y)}^{n+1}}
\begin{pmatrix}
\Gamma{(R+I)}\Gamma{(2R+I)^{-1}}& \lambda_1^{R}(x)  & \ldots & \lambda_{n-1}^{R}(x) & \lambda_n^{R}(x) \\
\gamma_0(y) & \gamma_1(y) & \ldots & \gamma_{n-1}(y) & \gamma_n(y) \\
0 & \gamma_0(y)  & \ldots & \binom{n-1}{n-2} \gamma_{n-2}(y) & \binom{n}{n-1} \gamma_{n-1}(y) \\
\vdots & \vdots  & \ddots & \vdots & \vdots \\
0 & 0 & \ldots & \gamma_0(y) & \binom{n}{1}\gamma_1(y)
\end{pmatrix}
\end{equation}
where,
\begin{equation}
\gamma_0(y) = \frac{1}{\phi_0(y)},
\end{equation}
\begin{equation}
\gamma_n(y) = -\frac{1}{\phi_0(y)} \left( \sum_{j=1}^{n} \binom{n}{j} \phi_j(y) \gamma_{n-j}(y) \right), \quad n = 1,2, \ldots
\end{equation}
\\
\textbf{Proof.} Consider two polynomial sequences $(\phi_n(y))_{n \in \mathbb{N}}$ and $ (\gamma_n(y))_{n \in \mathbb{N}}$ defined as:
\begin{equation}
\phi(y,t) = \phi_0(y) + \frac{t}{1!} \phi_1(y) + \frac{t^2}{2!} \phi_2(y) + \ldots + \frac{t^n}{n!} \phi_n(y) + \ldots, \quad n = 0,1, \ldots; \phi_0(y) \neq 0,
\end{equation}
\begin{equation}
\hat{\phi}(y,t) = \gamma_0(y) + \frac{t}{1!} \gamma_1(y) + \frac{t^2}{2!} \gamma_2(y) + \ldots + \frac{t^n}{n!} \gamma_n(y) + \ldots, \quad n = 0,1, \ldots; \gamma_0(y) \neq 0,
\end{equation}
satisfying
\begin{equation}
\phi(y,t) \hat{\phi}(y,t) = 1.
\end{equation}
Using the Cauchy-product rule, we find
\begin{equation}
\phi(y,t) \hat{\phi}(y,t) = \sum_{n=0}^{\infty} \sum_{j=0}^{n} \binom{n}{j} \phi_j(y) \gamma_{n-j}(y) \frac{t^n}{n!}.
\end{equation}
In view of (3.13) and (3.14), we get
\begin{equation}
\sum_{j=0}^{n} \binom{n}{j} \phi_j(y) \gamma_{n-j}(y) = \begin{cases} 
1 & \text{for } n = 0 \\
0 & \text{for } n > 0
\end{cases}.
\end{equation}
Therefore,
\[
\gamma_0(y) = \frac{1}{\phi_0(y)},
\]
\[
\gamma_n(y) = -\frac{1}{\phi_0(y)} \left( \sum_{j=1}^{n} \binom{n}{j} \phi_j(y) \gamma_{n-j}(y) \right), \quad n = 1,2, \ldots
\]
$\hat{\phi}(y,t)$ is multiplied on both sides of equation (2.6) and identity (3.13) is used to obtain
\begin{equation}
\hat\phi(y,t) \sum_{n=0}^{\infty} {_{G}\lambda_n^{R}(x,y)} \frac{t^n}{n!} = e^{xt} \cos{(\sqrt{t}; R)}.
\end{equation}
Using equations (1.11) and (3.12) in the previous equation, we obtain the following series by comparing the coefficients of comparable powers of $t$ on both sides of the resultant equation:
\begin{equation}
\sum_{j=0}^{n} \binom{n}{j} {_{G}\lambda_{n-j}^{R}(x,y)} \gamma_{j}(y) = \lambda_n^{R}(x),
\end{equation}
Again, using equation (1.6) gives the system of infinite equations with unknowns $_G\lambda_n^{R}(x,y), 
\\
\ (n = 0,1, \ldots)$:
\begin{equation}
\begin{aligned}
\begin{cases}
    &_G\lambda_0^{R}(x,y) \gamma_0(y) = \Gamma(R+1) \Gamma(2R + I)^{-1}, \\
    &_G\lambda_0^{R}(x,y) \gamma_1(y) + {_G\lambda_1}^{R}(x,y) \gamma_0(y)= \lambda_1^{R}(x), \\
    &_G\lambda_0^{R}(x,y) \gamma_2(y) + \binom{2}{1} {_G\lambda_1}^{R}(x,y) \gamma_1(y) +  {_G\lambda_2}^{R}(x,y) \gamma_0(y) = \lambda_2^{R}(x), \\
    &\vdots \\
     &_G\lambda_0^{R}(x,y) \gamma_{n-1}(y) +\binom{n-1}{n-2} {_G\lambda_{1}}^{R}(x,y) \gamma_{n-2}(y) + \cdots +{_G\lambda_{n-1}}^R(x,y) \gamma_0(y) = \lambda_{n-1}^{R}(x), \\
    &_G\lambda_0^{R}(x,y) \gamma_n(y) +\binom{n}{n-1} {_G\lambda_{1}}^{R}(x,y) \gamma_{n-1}(y) + \cdots +{_G\lambda_{n}}^R(x,y) \gamma_0(y) = \lambda_n^{R}(x). \\
    &\vdots \\
    \end{cases}
\end{aligned}
\end{equation}
The lower triangular form of the system provides a way to find the unknown $_G\lambda_n^R(x,y)$. Cramer's rule is used to evaluate the determinant in the denominator of the resultant equation and it is then transposed in the numerator as follows:
\begin{equation}
_G\lambda_n^{R}(x,y) = \frac{1}{\gamma_0{(y)}^{n+1}}
\begin{pmatrix}
\gamma_0(y) & \gamma_1(y)  & \ldots & \gamma_{n-1}(y) & \gamma_n(y) \\
0 & \gamma_0(y)  & \ldots & \binom{n-1}{n-2} \gamma_{n-2}(y) & \binom{n}{n-1} \gamma_{n-1}(y) \\
\vdots & \vdots & \ddots & \vdots & \vdots \\
0 & 0 & \ldots  & \gamma_0(y) & \binom{n}{1}\gamma_1(y)  \\
\Gamma{(R+I)}\Gamma{(2R+I)^{-1}}& \lambda_1^{R}(x) & \ldots & \lambda_{n-1}^{R}(x) & \lambda_n^{R}(x) \\
\end{pmatrix}
\end{equation}
which after performing n circular row exchanges, gives assertion (3.8).
 
 \vspace{+4pt} \noindent
 In next section, we discuss some members of the family of 2-variable general-$\lambda$-matrix polynomials 2VG$\lambda$MP and their properties. For different values of n and R, 3D graphs and zeros distribution of these polynomials are also shown.
\vspace{-12pt}
\section{Examples}
\vspace{-7pt}

The family of 2-variable general-$\lambda$-matrix polynomials 2VG$\lambda$MP ${_{G}\lambda_n^{R}(x,y)}$ encompasses a wide range of intriguing polynomials that exhibit unique mathematical properties. In our exploration, we investigate several notable members of this family by applying various symbolic operators, such as $\hat{q_y}$, to uncover their structural characteristics. Additionally,  a comprehensive analysis of specific instances of the family of 2-variable general-$\lambda$-matrix polynomials 2VG$\lambda$MP ${_{G}\lambda_n^{R}(x,y)}$ are performed, focusing on their various properties as well as 3D graphs and the distribution of their zeros. By examining these aspects, a deeper understanding of the complexities inherent in this fascinating class of polynomials are explored.

\pagebreak

\vspace{-10pt}
\subsection{Gould Hopper-$\lambda$-matrix polynomials ${_{H^m}\lambda_n^R(x,y)}$} 
\vspace{-4pt}
Gould Hopper polynomials ${H_n^m(x,y)}$ are defined in terms of the nth power of the binomial \cite{babusci2015lacunary}:
\begin{equation}
{H_n^m(x,y)}=(x+{_m\hat{b}_y})^n\phi_0,
\end{equation}
where, ${_m\hat{b}_y^r}$ is a symbolic operator defined on vacuum $\phi_0$\cite{babusci2015lacunary} such that
\begin{equation}
_m\hat{b}_y^r\phi_0=\frac{(y^{\frac{r}{m}})r!} {\Gamma{({\frac{r}{m}}+1)}}A_{m,r},
\end{equation}
\begin{equation}
A_{m,r} = \begin{cases} 
1, & r = mp,p \in \mathbb{N}, \\ 
0, & \text{otherwise}.
\end{cases}
\end{equation}
In view of the above operators, we derive some properties and results related to Gould Hopper-$\lambda$-matrix polynomials  ${_{H^m}\lambda_n^R(x,y)}$ which are as follows:

\begin{table}[h]

\small
\setlength{\tabcolsep}{4pt}  
\renewcommand{\arraystretch}{0.9} 
\begin{tabular*}{\textwidth}{@{} l l l @{}}
\hline
\textbf{No.} & \textbf{Result} & \textbf{Mathematical expression} \\
\hline
I & Umbral definition & \( _{H^m}\lambda_n^{R}(x,y)=\hat{J}^R(x+{_{m}\hat{b}_y}-\hat{J})^n\phi_0\psi_0 \) \\
II & Generating function & \( \sum_{n=0}^{\infty}{_{H^m}\lambda_n^{R}(x,y)}\frac{t^n}{n!}=e^{xt+yt^m}\cos({\sqrt{t},R}) \) \\
III & Series definition & \({_{H^m}\lambda_n^{R}(x,y)}=\sum_{j=0}^{n} \binom{n}{j} (-1)^jH_{n-j}^{(m)}(x,y)\frac{\Gamma{(R+(j+1)I)}}{\Gamma{(2R+(2j+1)I)}} \) \\
IV & Summation formulae & \(_{H^m}\lambda_n^{R}(x+z,y)=\sum_{j=0}^{n}\binom{n}{j}z^j{_{H^m}\lambda_{n-j}^{R}(x,y)}  \) \\
   &  & \(_{H^m}\lambda_{j+l}^{R}(z,y)=\sum_{n,r=0}^{j,l}\binom{j}{n}\binom{l}{r}(z-x)^{n+r}{_{H^m}\lambda_{j+l-n-r}^{R}(x,y)}\)\\
V & Multiplicative operator and & \((x+{_m\hat{b}_y}-\hat{J}) \) \\
  & Derivative Operator & \(D_x\)\\
VI & Symbolic recurrence relation & \( \frac{\partial}{\partial x}{_{H^m}\lambda_n^{R}(x,y)}=n{_{H^m}\lambda_{n-1}^{R}(x,y)} \) \\
   &  & \( \frac{\partial^j}{{\partial x}^j}{_{H^m}\lambda_n^{R}(x,y)}=\frac{n!}{(n-j)!} {_{H^m}\lambda_{n-j}^{R}(x,y)} \)\\
VII & Symbolic differential equation & \( ((x+{_m\hat{b}_y}-\hat{J})\frac{\partial}{\partial x}-n){_{H^m}\lambda_n^{R}(x,y)}=0 \) \\
\hline
\end{tabular*}

\caption{\textbf{Results for the Gould Hopper-$\lambda$-matrix polynomials}}

\end{table}

\noindent
Next, we conduct a graphical analysis of the Gould-Hopper-$\lambda$-matrix polynomials, focusing on the three-dimensional structure of the polynomial and the distribution of its zeros. This study involves examining various types of zeros, including real zeros and their specific patterns of arrangement. The analysis is illustrated through several figures that provide crucial insights into the behavior of the polynomial and its zeros:
Figure 1 presents a detailed 3D visualization of the polynomial, offering an overview of its structural properties; Figure 2 explores the distribution of the zeros, illustrating their precise locations in the complex plane; Figure 3 focuses on the real zeros, providing a detailed examination of their distribution along the real axis; Figure 4 presents a three-dimensional stacking representation of the zeros, highlighting both their real and imaginary components along with their corresponding index values.
This analysis provides valuable insights into the structure and behavior of Gould-Hopper $\lambda$-matrix polynomials, enhancing both theoretical understanding and potential applications in mathematical modeling and problem-solving.

\includegraphics[scale=0.615]{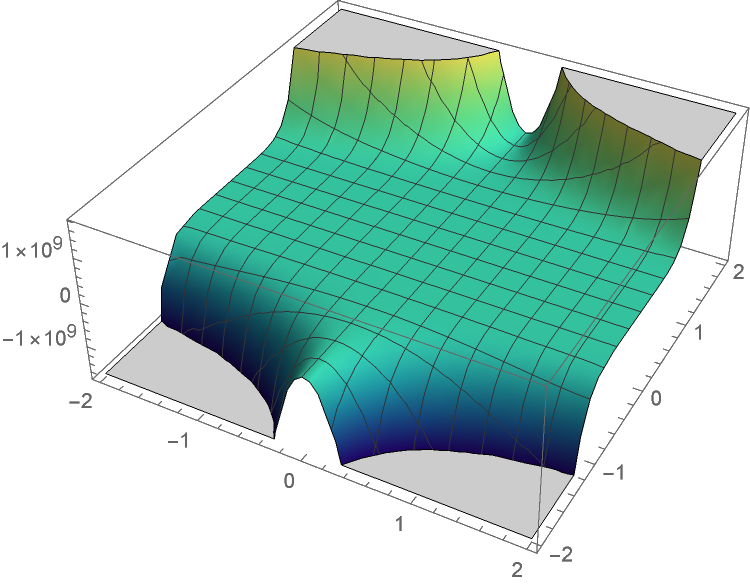}
\includegraphics[scale=0.615]{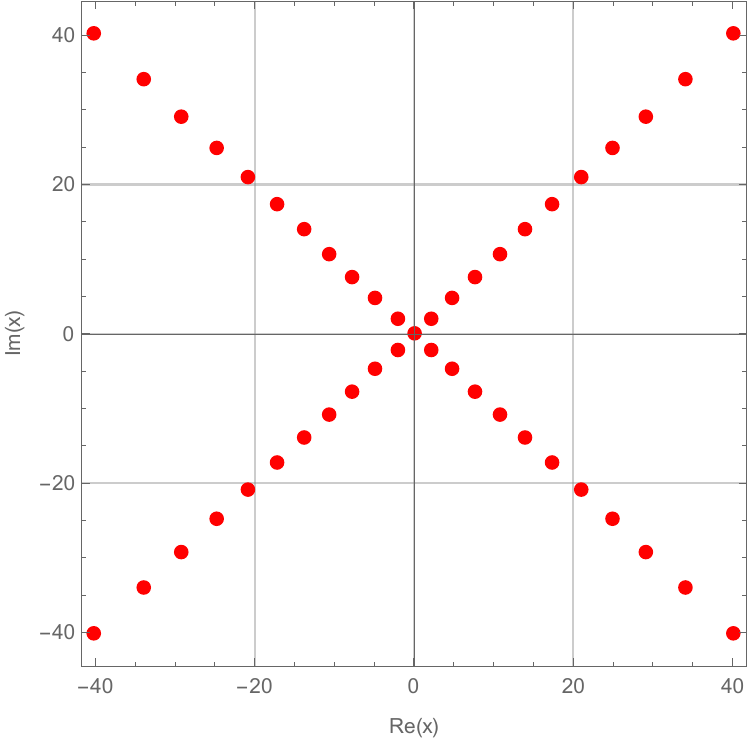}

\textbf{Figure 1.} 3D graph of ${_{H^4}\lambda_{30}^{15}(x,y)}$  \hspace{2.6cm} \textbf{Figure 2.} Zeros distribution of ${_{H^4}\lambda_{45}^{20}(x,6)}$ 

\includegraphics[scale=0.615]{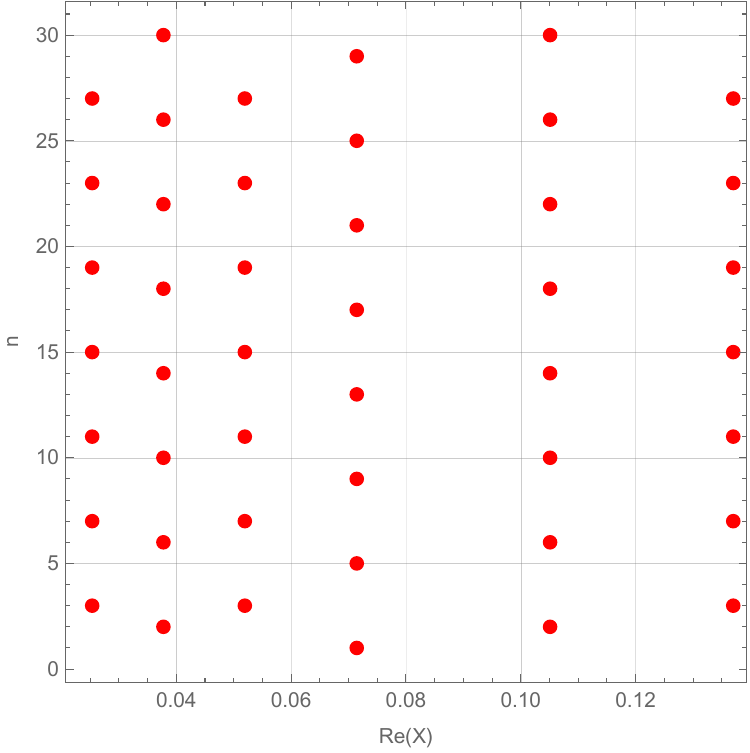}
\includegraphics[scale=0.615]{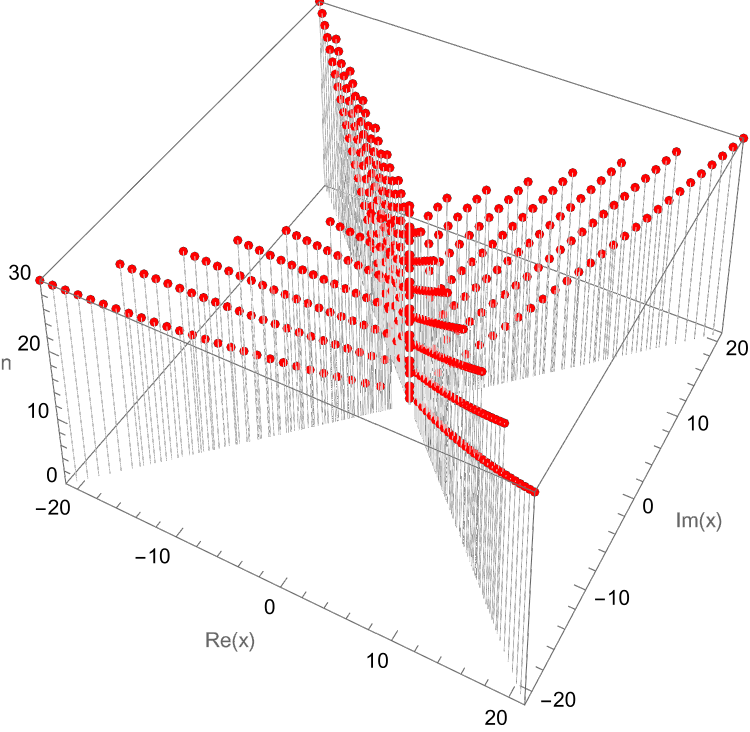}

\textbf{Figure 3.} Real zeros  of ${_{H^4}\lambda_{n}^{3}(x,1)}$ \hspace{2.6cm} \textbf{Figure 4.} Stacks of zeros of ${_{H^4}\lambda_{n}^{4}(x,2)}$ 

\vspace{+4pt} \noindent
\textbf{Remark 3.} Gould Hopper-$\lambda$-matrix polynomials for $m=2$ reduce to Hermite-$\lambda$-matrix polynomials and the results of Table 2 can be obtained analogously for Hermite-$\lambda$-matrix polynomials\cite{alatawi2024exploring} ${_{H}\lambda_n^{R}(x,y)}$.

\vspace{-10pt} 
\pagebreak
\subsection{Generalized Laguerre-$\lambda$-matrix polynomials ${_{mL}\lambda_n^R(x,y)}$}
\vspace{-4pt}
The umbral form of the generalized Laguerre polynomial ${_{m}{L}_n(y,x)}$ serves as a fundamental concept in the exploration of Generalized Laguerre-$\lambda$-matrix polynomials ${_{mL}\lambda_n^R(x,y)}$. This foundational framework provides a basis for further analysis and study of the Generalized Laguerre-$\lambda$-matrix polynomials. By utilizing the umbral form, we can delve into the properties, generating functions, explicit series representations and various identities associated with these matrix polynomials. First we define the Generalized Laguerre polynomials ${_{m}{L}_n(y,x)}$ by means of the following umbral definition:
\begin{equation}
{_{m}{L}_n(y,x)}=(x+{_{m}\hat{b}_{D_{y^{-1}}}})^n\phi_0,
\end{equation}
where, ${_{m}\hat{b}_{D_{y^{-1}}}^r}$ is a symbolic operator defined on vacuum $\phi_0$ such that 
\begin{equation}
{_{m}\hat{b}_{D_{y^{-1}}}^r}\phi_0=\frac{(D_{y^{-1}}^{\frac{r}{m}})r!} {\Gamma{({\frac{r}{m}}+1)}}A_{m,r},
\end{equation}
\begin{equation}
A_{m,r} = \begin{cases} 
1, & r = mp,p \in \mathbb{N}, \\ 
0, & \text{otherwise}.
\end{cases}
\end{equation}
Now, using the above relations, we find some properties and results related to Generalized Laguerre-$\lambda$-matrix polynomials ${_{mL}\lambda_n^R(x,y)}$ which are as follows:

\begin{table}[h]

\small
\setlength{\tabcolsep}{4pt}  
\renewcommand{\arraystretch}{0.9} 
\begin{tabular*}{\textwidth}{@{} l l l @{}}
\hline
\textbf{No.} & \textbf{Result} & \textbf{Mathematical expression} \\
\hline
I & Umbral definition & \( _{mL}\lambda_n^{R}(y,x)=\hat{J}^R(x+{_m{\hat{b}}_{D_{{y}^{-1}}}}-\hat{J})^n\phi_0\psi_0      \) \\
II & Generating function & \( \sum_{n=0}^{\infty}{_{mL}\lambda_n^{R}(y,x)}\displaystyle\frac{t^n}{n!}=e^{yt}C_0{(-yt^m)}\cos({\sqrt{t},R}) \) \\
III & Series definition & \({_{mL}\lambda_n^{R}(y,x)}=\sum_{j=0}^{n} \binom{n}{j} (-1)^j {_mL_{n-j}(y,x)}\frac{\Gamma{(R+(j+1)I)}}{(\Gamma{(2R+(2j+1)I))}} \) \\
IV & Summation formulae. & \(_{mL}\lambda_n^{R}(y+z,x)=\sum_{j=0}^{n}\binom{n}{j}z^j{_{mL}\lambda_{n-j}^{R}(y,x)}  \)  \\
 &  &\(_{mL}\lambda_{j+l}^{R}(z,x)=\sum_{n,r=0}^{j,l}\binom{j}{n}\binom{l}{r}(z-y)^{n+r}{_{mL}\lambda_{j+l-n-r}^{R}(y,x)}\)\\
V & Multiplicative operator & \((x+{_m{\hat{b}}_{D_{{y}^{-1}}}}-\hat{J}) \) \\
 & Derivative operator & \(D_x\)\\
VI & Symbolic recurrence relation & \( \displaystyle\frac{\partial}{\partial x}{_{mL}\lambda_n^{R}(y,x)}={(n)}{_{mL}\lambda_{n-1}^{R}(y,x)}    \)   \\ 
& & \(  \displaystyle\frac{\partial^j}{{\partial x}^j}{_{mL}\lambda_n^{R}(y,x)}=\displaystyle\frac{n!}{(n-j)!} {_{mL}\lambda_{n-j}^{R}(y,x)}   \)\\ 
VII & Symbolic differential equation & \(  ((x+{_m\hat{b}_{D_{{y}^{-1}}}}-\hat{J})\displaystyle\frac{\partial}{\partial x}-n){_{mL}\lambda_n^{R}(y,x)}=0 \) \\ 
\hline
\end{tabular*}
\caption{\textbf{Results for the Generalized Laguerre-$\lambda$-matrix polynomials}}
\end{table}
\noindent
Now, we examine the intricate details of the 3D graph, the distribution of zeros, real zeros and the distribution of zeros of the Generalized Laguerre-$\lambda$-matrix polynomials. These figures provide four key visualizations to illustrate these aspects:
Figure 5 presents a high-resolution surface plot that provides an in-depth 3D graphical representation of the structure of the polynomials. Figure 6 illustrates the distribution pattern of the zeros, clearly depicting their locations in the complex plane. Figure 7 focuses on the real zeros, plotting their distribution and behavior across various values of 'n' to reveal their underlying patterns. Lastly, Figure 8 offers a detailed visualization of the zeros, emphasizing their arrangement in the complex plane, along with their corresponding index values.
These visualizations are crucial for understanding the fundamental properties of the Generalized Laguerre-$\lambda$-matrix polynomials. By closely analyzing the arrangement and distribution of their zeros, we deepen our understanding of the behavior of this polynomial and their relevance in various mathematical contexts.

\vspace{+4pt}

\includegraphics[scale=0.615]{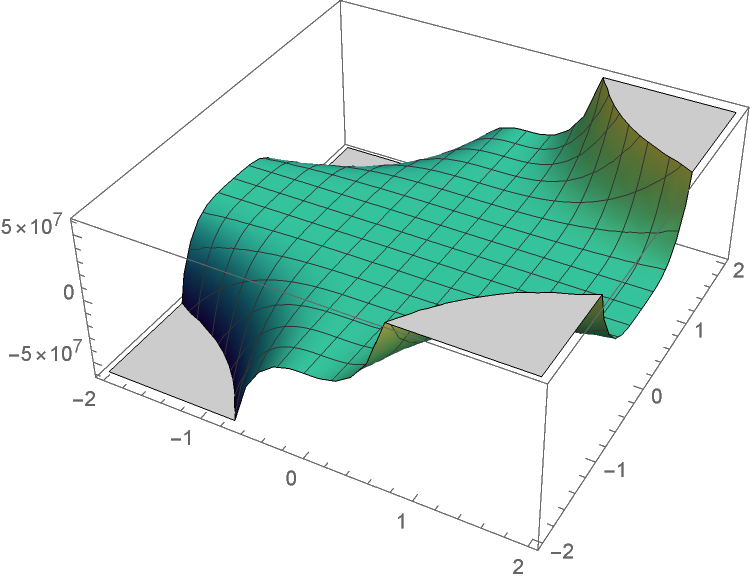}
\includegraphics[scale=0.615]{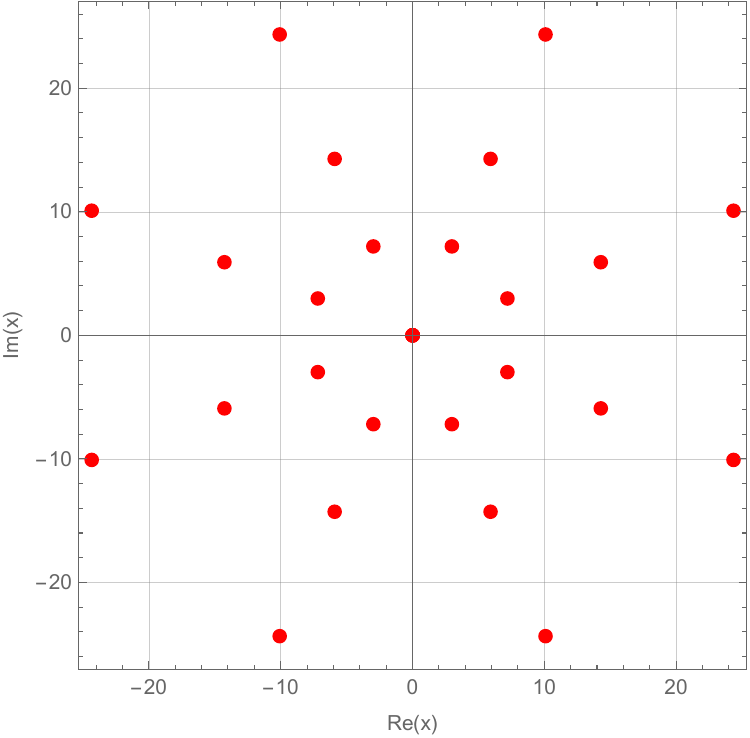}

\textbf{Figure 5.} 3D graph of ${_{8L}\lambda_{35}^{20}(x,y)}$ \hspace{2.6cm} \textbf{Figure 6.}  Zeros distribution of ${_{8L}\lambda_{30}^{25}(x,1)}$

\includegraphics[scale=0.615]{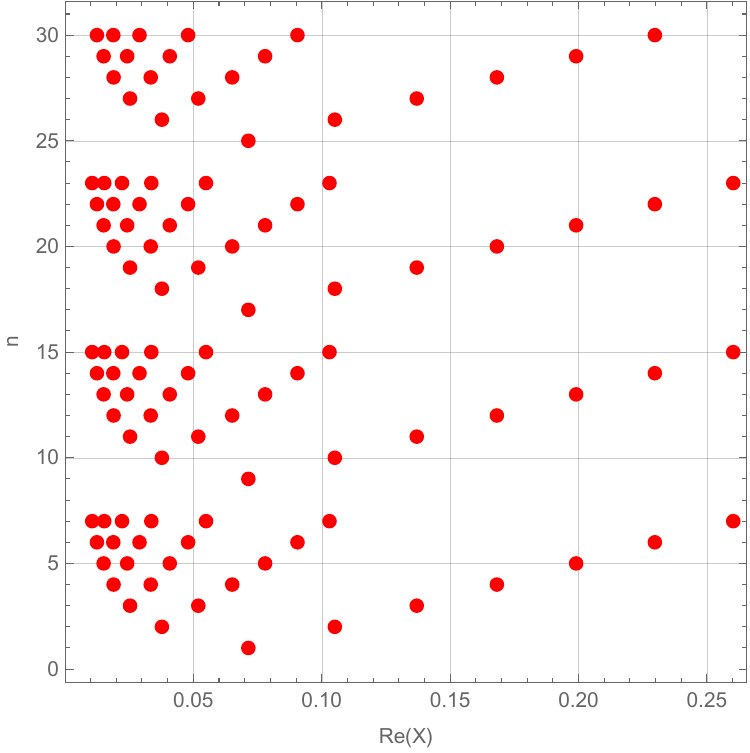}
\includegraphics[scale=0.615]{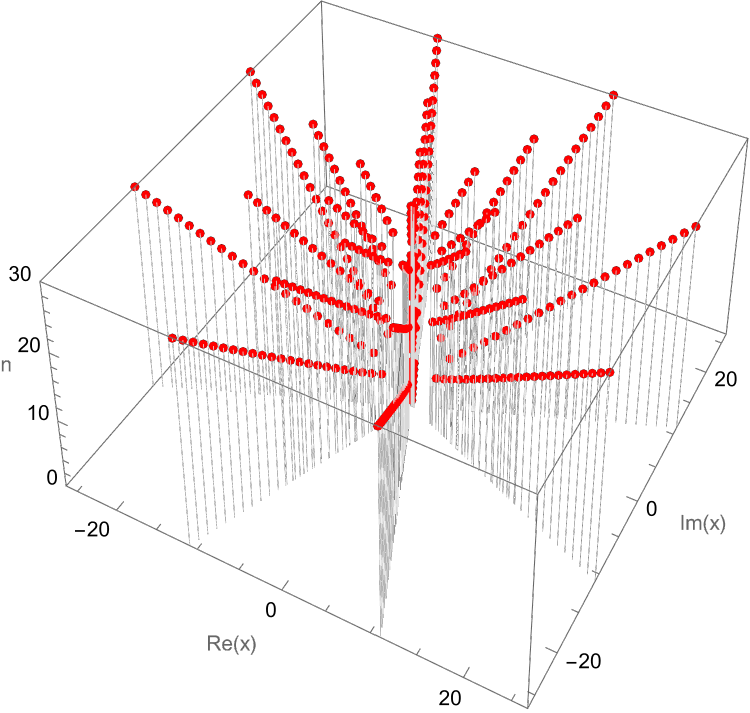}

\textbf{Figure 7.} Real zeros of ${_{8L}\lambda_{n}^{3}(x,1)}$ \hspace{2.6cm} \textbf{Figure 8.} Stacks of zeros of ${_{8L}\lambda_{n}^{2}(x,2)}$

\vspace{-10pt}
\pagebreak
\subsection{Truncated Exponential-$\lambda$-matrix polynomials ${_{e}\lambda_n^R(x,y)}$}
\vspace{-4pt}  
The umbral definition of  Truncated Exponential polynomials ${{e}_n(x,y)}$ is given by \cite{hibapreprint}
\begin{equation}
e_n(x,y)=(x+\hat{f}y)^n\beta_0,
\end{equation}
where, $\hat{f}$ is the umbral operator which operates on the vacuum $\beta$ such that 
\begin{equation}
\hat{f}^n\beta_0=\Gamma{(n+1)}.
\end{equation}
In case of Truncated Exponential-$\lambda$-Matrix polynomials, $\phi(y,t)=\frac{1}{1-yt}$ and $\phi_n(y)=n!y^n$

\vspace{+4pt} \noindent
In view of the above facts, we derive some properties of  Truncated Exponential-$\lambda$-Matrix polynomials ${_{e}\lambda_n^R(x,y)}$ which are as follows:

\begin{table}[h]

\small
\setlength{\tabcolsep}{4pt}  
\renewcommand{\arraystretch}{0.9} 
\begin{tabular*}{\textwidth}{@{} l l l @{}}
\hline
\textbf{No.} & \textbf{Result} & \textbf{Mathematical expression} \\
\hline
I & Umbral definition & \( _{e}\lambda_n^{R}(x,y)=\hat{J}^R(x+\hat{f}y-\hat{J})^n\beta_0\psi_0      \) \\
II & Generating function & \( \sum_{n=0}^{\infty}{_{e}\lambda_n^{R}(x,y)}\displaystyle\frac{t^n}{n!}=e^{xt}\frac{1}{1-yt}\cos({\sqrt{t},R}) \) \\
III & Series definition & \({_{e}\lambda_n^{R}(x,y)}=\sum_{j=0}^{n} \binom{n}{j} (-1)^j e_{n-j}(x,y)\frac{\Gamma{(R+(j+1)I)}}{(\Gamma{(2R+(2j+1)I))}} \) \\
IV & Summation formulae & \(_{e}\lambda_n^{R}(x+z,y)=\sum_{j=0}^{n}\binom{n}{j}z^j{_{e}\lambda_{n-j}^{R}(x,y)}\) \\
 &  &\(_{e}\lambda_{j+l}^{R}(z,y)=\sum_{n,r=0}^{j,l}\binom{j}{n}\binom{l}{r}(z-x)^{n+r}{_{e}\lambda_{j+l-n-r}^{R}(x,y)}\)\\
V & Multiplicative operator and & \((x+\hat{f}y-\hat{J}) \) \\
 & Derivative operator & \(D_x\) \\
VI & Symbolic recurrence relation & \( \displaystyle\frac{\partial}{\partial x}{_{e}\lambda_n^{R}(x,y)}=n{_{e}\lambda_{n-1}^{R}(x,y)}    \)  \\ 
& & \(  \displaystyle\frac{\partial^j}{{\partial x}^j}{_{e}\lambda_n^{R}(x,y)}=\displaystyle\frac{n!}{(n-j)!} {_{e}\lambda_{n-j}^{R}(x,y)}   \) \\
VII & Symbolic differential equation & \(  ((x+\hat{f}y-\hat{J})\displaystyle\frac{\partial}{\partial x}-n){_{e}\lambda_n^{R}(x,y)}=0 \) \\ 
\hline
\end{tabular*}
\caption{\textbf{Results for the Truncated Exponential-$\lambda$-matrix polynomials}}
\end{table}
\noindent
Next, we'll explore the 3D graphs,the distribution of zeros, real zeros and stacks of zeros for Truncated Exponential-$\lambda$-matrix polynomials. Figures 9, 10, 11, and 12 provide detailed representations of these concepts.
Figure 9 showcases the 3D graph, vividly illustrating the polynomial's behavior across different dimensions, revealing its overall structure. This visualization helps us understand how the polynomial behaves under various conditions.
Meanwhile, Figure 10 provides insight into the distribution of zeros, highlighting where the polynomial intersects the axis. Analyzing the zeros is crucial, as they provide valuable information about the roots of the polynomial and their implications in various applications. 
Figure 11 shows Real zeros of ${_{e}\lambda_{n}^{R}(x,y)} $.
Figure 12 displays the stacks of zeros. 
Together, these figures offer a comprehensive view of the mathematical characteristics of Truncated Exponential-$\lambda$-matrix polynomials, paving the way for deeper analysis and interpretation.

\vspace{+4pt} \noindent

\includegraphics[scale=0.615]{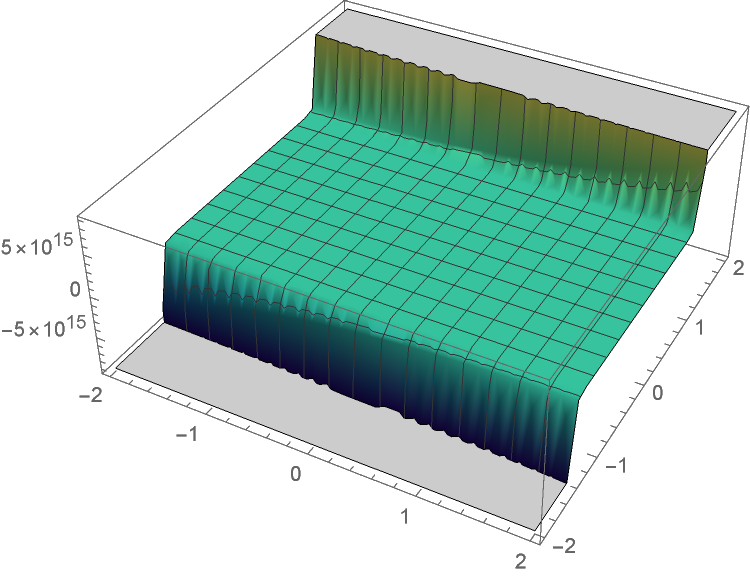}
\includegraphics[scale=0.615]{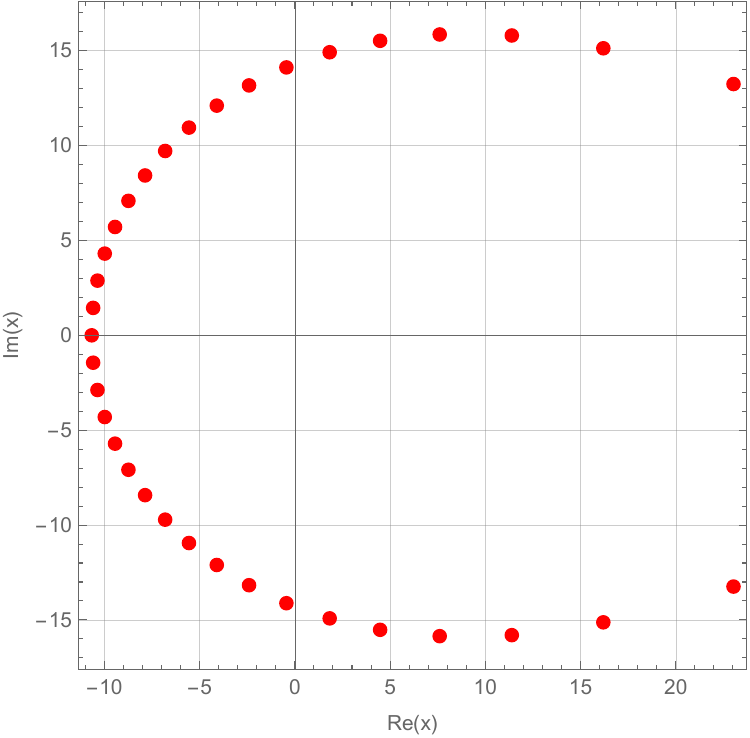}

\textbf{Figure 9.} 3D graph of ${_{e}\lambda_{31}^{16}(x,y)} $ \hspace{2.6cm} \textbf{Figure 10.} Zeros distribution of ${_{e}\lambda_{35}^{10}(x,1)}$

\includegraphics[scale=0.615]{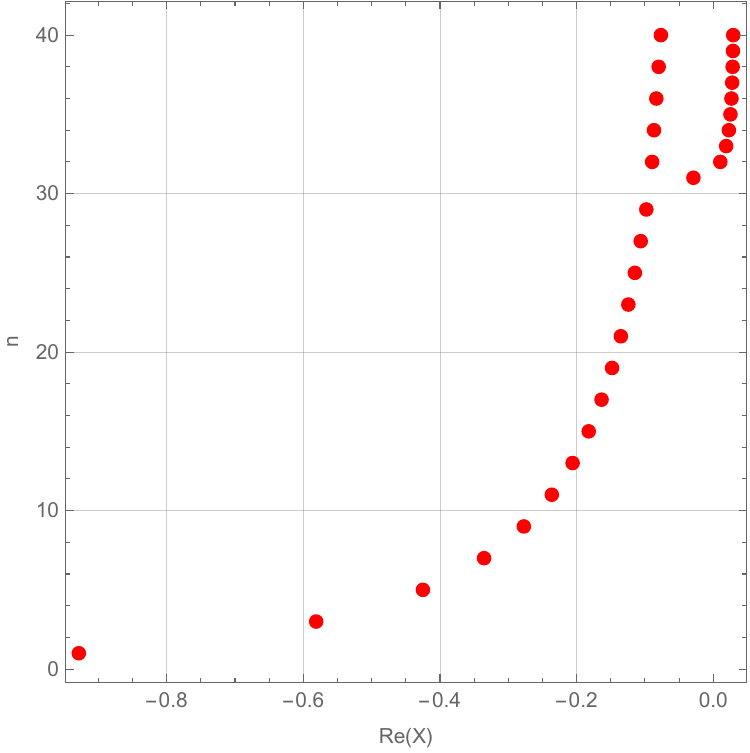}
\includegraphics[scale=0.615]{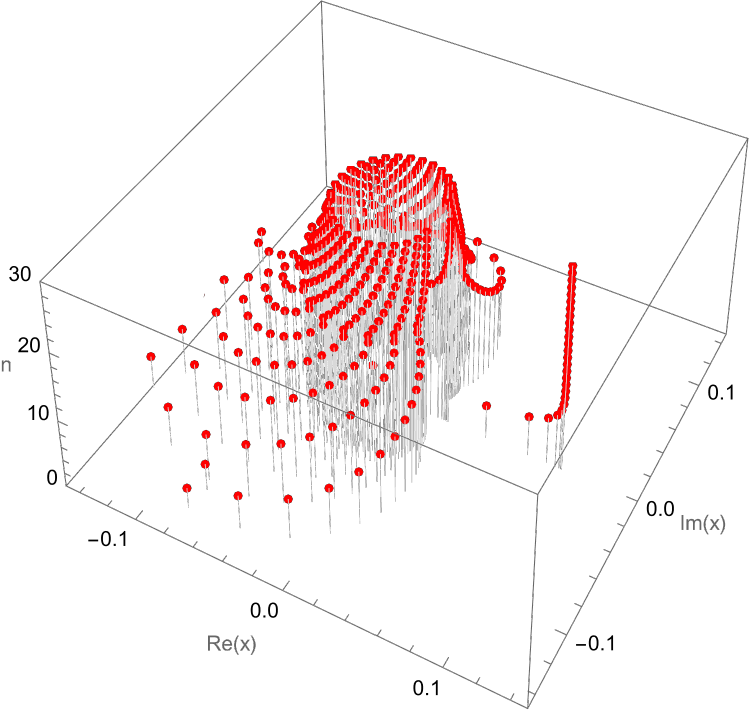}

\textbf{Figure 11.} Real zeros of ${_{e}\lambda_{n}^{3}(x,1)}$ \hspace{2.6cm} \textbf{Figure 12.} Stacks of zeros of ${_{e}\lambda_{n}^{1}(x,1/2)}$

\vspace{-12pt}
\pagebreak
\section{Conclusion}
\vspace{-7pt}
The potential applications of 2-variable general-$\lambda$-matrix polynomials for future research are diverse. This paper concludes by exploring the new extension of these polynomials as q-2-variable general-$\lambda$-matrix polynomials (q2VG$\lambda$MP) ${_{G}\lambda_{n,q}^{R}(x,y)}$. There are numerous avenues for conducting a detailed study on q-2-variable general-$\lambda$-matrix polynomials. In this section, we focus on examining one member of the 2-variable general polynomials, i.e., the 2-variable Hermite polynomials. The conclusion of this work involves an examination of the q-Hermite $\lambda$-matrix polynomials (qH$\lambda$MP) ${_{H}\lambda_{n,q}^{R}(x,y)}$ as a member of the q-2-variable general $\lambda$-matrix polynomials family.

\vspace{+4pt} \noindent
The following is the definition of the q gamma function\cite{ismail1981basic}:
\begin{equation}
    \Gamma_q(x+1)={[x]_q}\Gamma_q(x),
\end{equation}
q analogue of $x$ is provided by \cite{andrews1999special}:
\begin{equation}
    {[x]_q}=\frac{1-q^x}{1-q}.
\end{equation}
The definition of the Nalli-Ward-Alsalam q-addition (NWA) is given below\cite{al1958q}:
\begin{equation}
(a \oplus_q b)^n=\sum_{r=0}^{n}{\binom{n}{r}_q} a^r b^{n-r}.
\end{equation}
Introducing the q-Hermite $\lambda$-matrix polynomials involves the introduction of certain operators and functions.
\\
Here, we consider a symbolic operator $\hat{J}_q$, which acts on vacuum $\psi_0$ such that:
\begin{equation}
    \hat{J}_q^R\psi_0=\frac{\Gamma_q{(R+(j+1)I)}}{\Gamma_q{(2R+(2j+1)I)}}
\end{equation}
and
\begin{equation}
    \hat{J}_q^{R+P}\psi_0=\hat{J}_q^R\hat{J}_q^P\psi_0,
\end{equation}
where R and P are positive stable matrices.

\vspace{+4pt} \noindent
We define the q-Hermite function vacuum as follows:
\begin{equation}
     \theta_q(r)=
     \theta_{q,r}=
     \Biggl\{
     \begin{aligned}
     && y^{r}\frac{{(2r)_q!}}{(r)_q!} ,&&& z=2r,
     \\
     && 0                               ,&&& z=2r+1,
     \end{aligned}  
\end{equation}
\\
Also, consider the umbral operator $_y\hat{h}_q^r$ which acts on the vacuum $\theta_0$ such that $\forall m \in \mathbb{R},$ 
\begin{equation}
_y\hat{h}_q^m\theta_0=\theta_{q,m}.
\end{equation}
\\
In view of the equations (5.3), (5.6) and (5.7), we establish the umbral form of q-Hermite polynomials as follows:
\begin{equation}
    H_{n,q}(x.y)=(_y\hat{h}_q \oplus x)^n\theta_0.
\end{equation}
Next, we introduce umbral form of q-Hermite $\lambda$-matrix polynomials qH$\lambda$MP ${_{H}\lambda_{n,q}^{R}(x,y)}$ as:

\begin{equation}
    {_H\lambda_{n,q}^R(x,y)}=\hat{J}_q^R[(-\hat{J}_q) \oplus_q x \oplus_q {_y\hat{h}_q}]^n\theta_0\psi_0.
\end{equation}
\\
Applying equation (5.3) and expanding the r.h.s. of equation (5.9), we obtain the following expression:
  \begin{equation}
  {_{H}\lambda_{n,q}^{R}(x,y)}=(\hat{J}_q)^R\sum_{j=0}^{n}{\binom{n}{j}}_q (-\hat{J}_q)^j {(x \oplus _y\hat{q}_q)}^{n-j}\theta_0\psi_0.
\end{equation}
Using equations (5.8), (5.4), and (5.5), we derive the series definition of the q-Hermite $\lambda$-matrix polynomials, which is defined as follows:
\begin{equation}
 {_{H}\lambda_{n,q}^{R}(x,y)}=\sum_{j=0}^{n} {\binom{n}{j}}_q (-1)^j H_{n-j,q}(x,y)\Gamma_q{(R+(j+1)I)}(\Gamma_q{(2R+(2j+1)I))^{-1}}.   
\end{equation}
\noindent
Next, the graphs illustrate the 3D surface plot, distribution of zeros, real zeros, and stacks of zeros for the q-Hermite $\lambda$-matrix polynomials, as shown in Figures 13, 14, 15, and 16.
This visualization is essential for analyzing and comprehending the behavior and characteristics of these polynomials, represented as ${_H\lambda_{n,q}^R(x,y)}$.
The 3D graph will illustrate the complex relationships and patterns within the polynomial's structure, while the distribution of zeros will shed light on the locations of their roots in the complex plane. 

\vspace{+4pt} \noindent

\includegraphics[scale=0.615]{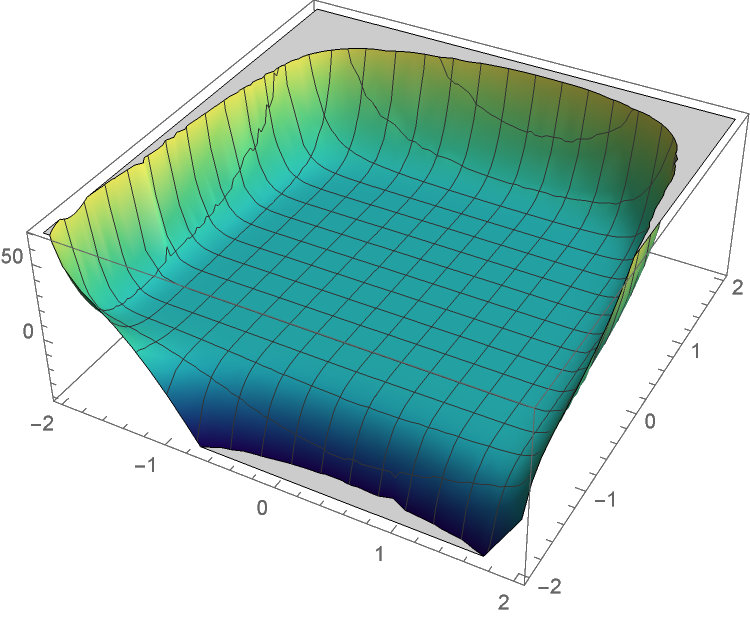}
\includegraphics[scale=0.615]{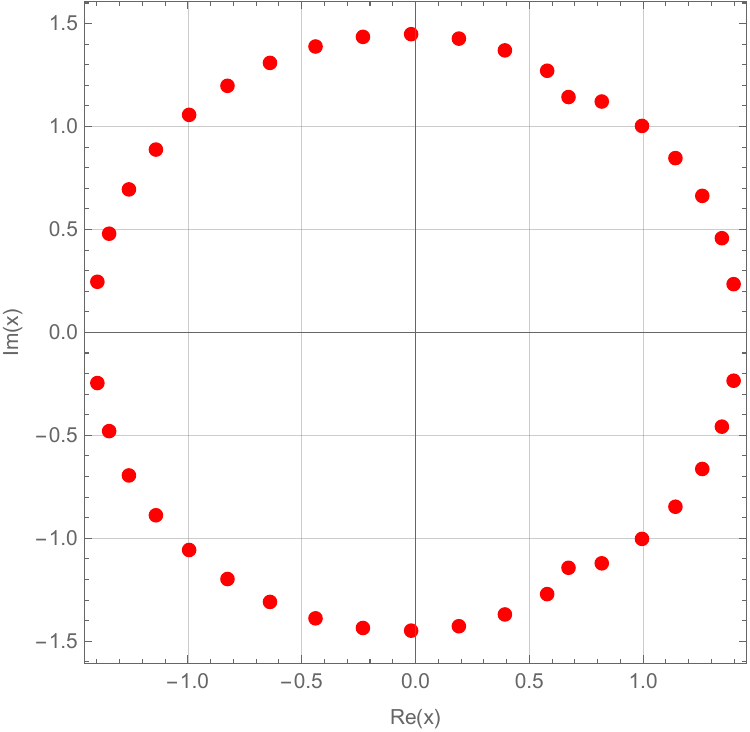}

\textbf{Figure 13.} D-Graph of ${_{H}\lambda_{22,1/2}^{15}(x,y)}$ \hspace{1cm} \textbf{Figure 14.} Zeros distribution of ${_{H}\lambda_{40,1/2}^{20}(x,1)}$

\includegraphics[scale=0.615]{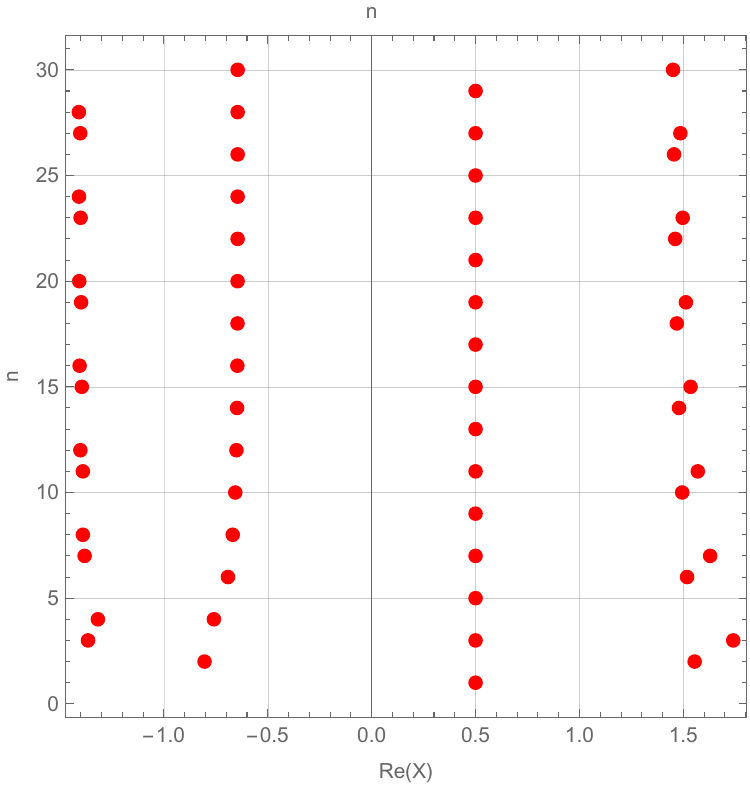}
\includegraphics[scale=0.615]{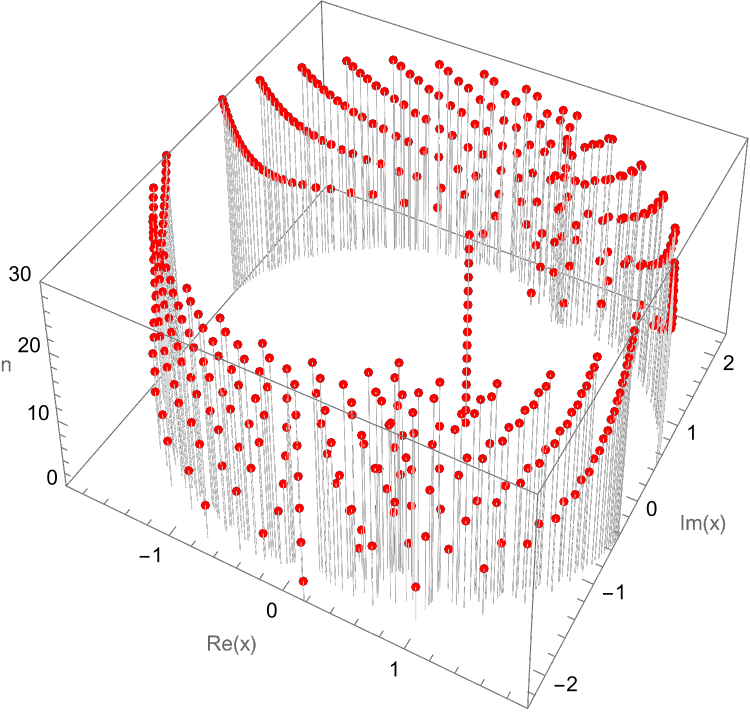}

\textbf{Figure 15.} Real zeros of ${_{H}\lambda_{n,1/2}^{20}(x,-1)}$ \hspace{1cm}\textbf{Figure 16.} Stacks of zeros of ${_{H}\lambda_{n,1/2}^{4}(x,2)}$

\vspace{+4pt} \noindent
By exploring these elements, we can better understand how the q-Hermite $\lambda$-matrix polynomials function under varying conditions and parameters, emphasizing their distinctive properties and potential applications in both mathematical theory and practical scenarios. This detailed representation will be a foundation for further investigation and research in this field.

\vspace{-12pt}
\bibliographystyle{plain}
\bibliography{mybib.bib}

\begin{thebibliography}{10}

\bibitem{al1958q}
Waleed~A Al-Salam.
\newblock q-bernoulli numbers and polynomials.
\newblock {\em Mathematische Nachrichten}, 17(3-6):239--260, 1958.

\bibitem{alatawi2024exploring}
Maryam~Salem Alatawi, Manoj Kumar, Nusrat Raza, and Waseem~Ahmad Khan.
\newblock Exploring zeros of hermite-$\lambda$ matrix polynomials: A numerical approach.
\newblock {\em Mathematics}, 12(10):1497, 2024.

\bibitem{andrews1999special}
George~E Andrews, Richard Askey, Ranjan Roy, Ranjan Roy, and Richard Askey.
\newblock {\em Special functions}, volume~71.
\newblock Cambridge university press Cambridge, 1999.

\bibitem{appell1926fonctions}
Paul Appell and Joseph~Kamp{\'e} De~F{\'e}riet.
\newblock {\em Fonctions hyperg{\'e}om{\'e}triques et hypersph{\'e}riques: polynomes d'Hermite}.
\newblock Gauthier-villars, 1926.

\bibitem{babusci2015lacunary}
D~Babusci, G~Dattoli, Bruna Germano, Maria~Renata Martinelli, PE~Ricci, et~al.
\newblock Lacunary generating functions of hermite polynomials and symbolic methods.
\newblock {\em Ilirias Journal of Mathematics}, 4:16--23, 2015.

\bibitem{dattoli1999hermite}
G~Dattoli.
\newblock Hermite-bessel and laguerre-bessel functions: A by-product of the monomiality principle.
\newblock {\em Advanced special functions and applications}, 1:147--164, 1999.

\bibitem{dattoli1999generalized}
G~Dattoli, S~Lorenzutta, AM~Mancho, and A~Torre.
\newblock Generalized polynomials and associated operational identities.
\newblock {\em Journal of computational and applied mathematics}, 108(1-2):209--218, 1999.

\bibitem{dattoli2004class}
G~Dattoli, Mauro Migliorati, HM~Srivastava, et~al.
\newblock A class of bessel summation formulas and associated operational methods.
\newblock {\em Fractional Calculus and Applied Analysis}, 7(2):169--176, 2004.

\bibitem{dattoli1998operational}
G~Dattoli and A~Torre.
\newblock Operational methods and two variable laguerre polynomials.
\newblock {\em Atti Accad. Sci. Torino Cl. Sci. Fis. Mat. Natur}, 132:3--9, 1998.

\bibitem{dattoli2017circular}
Giuseppe Dattoli, Emanuele Di~Palma, Silvia Licciardi, and Elio Sabia.
\newblock From circular to bessel functions: a transition through the umbral method.
\newblock {\em Fractal and Fractional}, 1(1):9, 2017.

\bibitem{dunford1988linear}
Nelson Dunford and Jacob~T Schwartz.
\newblock {\em Linear operators, part 1: general theory}, volume~10.
\newblock John Wiley \& Sons, 1988.

\bibitem{gould1962operational}
Henry~W Gould and AT~Hopper.
\newblock Operational formulas connected with two generalizations of hermite polynomials.
\newblock 1962.

\bibitem{ismail1981basic}
Mourad~EH Ismail.
\newblock The basic bessel functions and polynomials.
\newblock {\em SIAM Journal on Mathematical Analysis}, 12(3):454--468, 1981.

\bibitem{khan2013general}
Subuhi Khan and Nusrat Raza.
\newblock General-appell polynomials within the context of monomiality principle.
\newblock {\em International Journal of Analysis}, 2013(1):328032, 2013.

\bibitem{louisell1973quantum}
William~Henry Louisell.
\newblock Quantum statistical properties of radiation.
\newblock 1973.

\bibitem{steffensen1941poweroid}
JF5953 Steffensen.
\newblock The poweroid, an extension of the mathematical notion of power.
\newblock 1941.

\bibitem{hibapreprint}
Ghazala Yasmin and Hibah Islahi.
\newblock Truncated polynomial: Operational versus umbral methods (preprint).

\bibitem{zainab2024symbolic}
U~Zainab and N~Raza.
\newblock The symbolic approach to study the family of appell-$\lambda$ matrix polynomials.
\newblock {\em Filomat}, 38(4):1291--1304, 2024.

\end{thebibliography}
\end{document}